# ROBUST RECONSTRUCTION ON TREES IS DETERMINED BY THE SECOND EIGENVALUE

By Svante Janson and Elchanan Mossel[1]

*Uppsala University and University of California, Berkeley*

Consider a Markov chain on an infinite tree $T = (V, E)$ rooted at $\rho$. In such a chain, once the initial root state $\sigma(\rho)$ is chosen, each vertex iteratively chooses its state from the one of its parent by an application of a Markov transition rule (and all such applications are independent). Let $\mu_j$ denote the resulting measure for $\sigma(\rho) = j$. The resulting measure $\mu_j$ is defined on configurations $\sigma = (\sigma(x))_{x \in V} \in \mathcal{A}^V$, where $\mathcal{A}$ is some finite set. Let $\mu_j^n$ denote the restriction of $\mu$ to the sigma-algebra generated by the variables $\sigma(x)$, where $x$ is at distance exactly $n$ from $\rho$. Letting $\alpha_n = \max_{i,j \in \mathcal{A}} d_{\mathrm{TV}}(\mu_i^n, \mu_j^n)$, where $d_{\mathrm{TV}}$ denotes total variation distance, we say that the *reconstruction problem is solvable* if $\liminf_{n \to \infty} \alpha_n > 0$. Reconstruction solvability roughly means that the $n$th level of the tree contains a nonvanishing amount of information on the root of the tree as $n \to \infty$.

In this paper we study the problem of *robust reconstruction*. Let $\nu$ be a nondegenerate distribution on $\mathcal{A}$ and $\varepsilon > 0$. Let $\sigma$ be chosen according to $\mu_j^n$ and $\sigma'$ be obtained from $\sigma$ by letting for each node independently, $\sigma(v) = \sigma'(v)$ with probability $1 - \varepsilon$ and $\sigma'(v)$ be an independent sample from $\nu$ otherwise. We denote by $\mu_j^n[\nu, \varepsilon]$ the resulting measure on $\sigma'$. The measure $\mu_j^n[\nu, \varepsilon]$ is a perturbation of the measure $\mu_j^n$. Letting $\alpha_n(\nu, \varepsilon) = \max_{i,j \in \mathcal{A}} d_{\mathrm{TV}}(\mu_i^n[\nu, \varepsilon], \mu_j^n[\nu, \varepsilon])$, we say that the reconstruction problem is $\nu$-*robust-solvable* if $\liminf_{n \to \infty} \alpha_n(\nu, \varepsilon) > 0$ for all $0 < \varepsilon < 1$. Roughly speaking, the reconstruction problem is robust-solvable if for any noise-rate and for all $n$, the $n$th level of the tree contains a nonvanishing amount of information on the root of the tree.

Standard techniques imply that if $T$ is the rooted $B$-ary tree (where each node has $B$ children) and if $B|\lambda_2(M)|^2 > 1$, where $\lambda_2(M)$ is the second largest eigenvalue of $M$ (in absolute value), then for all nondegenerate $\nu$, the reconstruction problem is $\nu$-robust-solvable. We

---

Received December 2002; revised December 2003.

[1]Supported by a Miller fellowship at University of California, Berkeley.

*AMS 2000 subject classifications.* Primary 60K35; secondary 60J80, 82B26.

*Key words and phrases.* Robust phase transition, reconstruction on trees, branching number.

---

prove a converse and show that the reconstruction problem is not $\nu$-robust-solvable if $B|\lambda_2(M)|^2 < 1$. This proves a conjecture by the second author and Y. Peres. We also consider other models of noise and general trees.

**1. Introduction.** In this paper we study the perturbative theory of reconstruction on trees, and show how it depends on the spectrum of the underlying Markov chain. In particular, we show that the threshold for "robust reconstruction" for the $B$-ary tree is $B|\lambda_2(M)|^2 = 1$, where $\lambda_2(M)$ denotes the eigenvalue of $M$ which is the second largest in absolute value. In Section 3 we prove a similar threshold for general bounded degree trees, where $B$ is replaced by the branching number of the tree $\mathrm{br}(T)$. We refer the reader to Section 1.2 and to [3, 7, 21, 22] for background.

1.1. *Definitions and main results.* We proceed with some formal definitions. Let $T = (V, E, \rho)$ be a tree $T$ with nodes $V$, edges $E$ and root $\rho \in V$. We direct all edges away from the root so that if $e = (x, y)$, then $x$ is on the path connecting $\rho$ to $y$. Let $d(\cdot, \cdot)$ denote the graph-metric distance on $T$, and $L_n = \{v \in V : d(\rho, v) = n\}$ be the $n$th level of the tree. For $x \in V$ and $e = (y, z) \in E$, we denote $|x| = d(\rho, x)$, $d(x, (y, z)) = \max\{d(x, y), d(x, z)\}$ and $|e| = d(\rho, e)$. The $B$-*ary* tree is the infinite rooted tree, where each vertex has exactly $B$ children.

A Markov chain on the tree is a probability measure whose state space is $\mathcal{A}^V$, where $\mathcal{A}$ is a finite set. Without loss of generality we assume that $\mathcal{A} = \{1, \ldots, q\}$. Assume first that $T$ is finite and let $M = (M_{i,j})_{i,j \in \mathcal{A}}$ be a stochastic matrix. In this case the probability measure defined by $M$ on $T$ is given by

$$(1) \qquad \bar{\mu}_\ell(\sigma) = \mathbf{1}_{\{\sigma(\rho) = \ell\}} \prod_{(x,y) \in E} M_{\sigma(x), \sigma(y)}.$$

In other words, in $\bar{\mu}_\ell$ the root state $\sigma(\rho)$ satisfies $\sigma(\rho) = \ell$ and then each vertex iteratively chooses its state from the one of its parent by an application of the Markov transition rule given by $M$ (and all such applications are independent). We can define the measure $\bar{\mu}_\ell$ on an infinite tree too, by Kolmogorov's extension theorem, but we will not need chains on infinite trees in this paper (see [7] for basic properties of Markov chains on trees).

Instead, for an infinite tree $T$, we let $T_n = (V_n, E_n, \rho)$, where $V_n = \{x \in V : d(x, \rho) \leq n\}$, $E_n = \{e \in E : d(e, \rho) \leq n\}$ and define $\bar{\mu}_\ell^n$ by (1) for $T_n$. More explicitly,

$$(2) \qquad \bar{\mu}_\ell^n(\sigma) = \mathbf{1}_{\{\sigma(\rho) = \ell\}} \prod_{(x,y) \in E_n} M_{\sigma(x), \sigma(y)}.$$

We are particularly interested in the distribution of the states $\sigma(x)$ for $x \in L_n$, the set of leaves in $T_n$. This distribution, denoted by $\mu_k^n$, is the projection of $\bar{\mu}_k^n$ on $\mathcal{A}^{L_n}$ given by

$$(3) \qquad \mu_k^n(\sigma) = \sum \{\bar{\mu}_k^n(\bar{\sigma}) : \bar{\sigma}|L_n = \sigma\}.$$

In this paper we are interested in perturbative theory of the above process. Below we give three definitions of perturbations of $\mu_k^n$ representing three different types of "noise." We call a distribution $\nu$ on $\mathcal{A} = \{1, \ldots, q\}$ nondegenerate, if $\nu(i) > 0$, for all $1 \le i \le q$.

In the general setting the perturbation is obtained by observing, for leaves $x \in L_n$, not the state $\sigma(x)$ but a state (in a state space $\mathcal{B}$ possibly different from $\mathcal{A}$) derived from $\sigma(x)$ by another random choice (independently for all leaves). The extra choice can be described by a stochastic matrix $N = (N_{i,j})_{i \in \mathcal{A}, j \in B}$; this defines a probability measure on $\mathcal{A}^{V_n} \times \mathcal{B}^{L_n}$ by

$$(4) \qquad \bar{\mu}[N]_\ell^n(\sigma, \tau) = \mathbf{1}_{\{\sigma(\rho) = \ell\}} \prod_{(x,y) \in E_n} M_{\sigma(x),\sigma(y)} \times \prod_{y \in L_n} N_{\sigma(y),\tau(y)},$$

and the distribution of our observed states is the projection $\mu[N]_\ell^n$ on $\mathcal{B}^{L_n}$ given by

$$(5) \qquad \mu[N]_\ell^n(\tau) = \sum_\sigma \bar{\mu}[N]_\ell^n(\sigma, \tau).$$

We will mostly be interested in the following types of noise:

- Given $k \ge 0$, define $N = M^k$. Here, for each leaf independently, $k$ additional steps of the chain are performed. We write $\mu_\ell^n[k]$ for $\mu_\ell^n[N]$.
- Given a distribution $\nu$ on $\mathcal{A}$, define $N_{i,j} = (1 - \varepsilon)\mathbf{1}_{\{i=j\}} + \varepsilon\nu_j$. Here, for each leaf independently, with probability $1 - \varepsilon$, there is no noise; otherwise, the leaf state is chosen independently from anything else according to $\nu$. We will write $\mu_k^n[\nu, \varepsilon]$ for $\mu[N]_\ell^n$.
- Given $0 \le \epsilon \le 1$, we let $N$ be a $q \times (q+1)$ matrix defined by $N_{i,i} = (1-\epsilon)$, $N_{i,q+1} = \epsilon$ and $N_{i,j} = 0$ otherwise. Here, for each leaf independently, the state at the leaf is deleted with probability $\epsilon$ (deletion is marked by $q+1$). We write $\mu_\ell^n[\epsilon]$ for $\mu_\ell^n[N]$.

Recall that for distributions $\mu$ and $\nu$ on the same space $\Omega$, the total variation distance between $\mu$ and $\nu$ is

$$(6) \qquad D_V(\mu, \nu) = \tfrac{1}{2} \sum_{\sigma \in \Omega} |\mu(\sigma) - \nu(\sigma)|.$$

DEFINITION 1.1. (i) The reconstruction problem for the $B$-ary tree $T$ and $M$ is *solvable* if there exist $i, j \in \mathcal{A}$, for which

$$(7) \qquad \liminf_{n \to \infty} D_V(\mu_i^n, \mu_j^n) > 0.$$

(ii) The reconstruction problem for the $B$-ary tree $T$ and $M$ is *robust-solvable* if for all $k < \infty$ there exist $i, j \in \mathcal{A}$ for which

$$(8) \qquad \liminf_{n \to \infty} D_V(\mu_i^n[k], \mu_j^n[k]) > 0.$$

(iii) Let $\nu$ be a nondegenerate distribution. The reconstruction problem for the $B$-ary tree $T$ and $M$ is *$\nu$-robust-solvable* if for all $\varepsilon < 1$, there exist $i, j \in \mathcal{A}$, for which

$$(9) \qquad \liminf_{n \to \infty} D_V(\mu_i^n[\nu, \varepsilon], \mu_j^n[\nu, \varepsilon]) > 0.$$

(iv) The reconstruction problem for the $B$-ary tree $T$ and $M$ is *erasure-robust-solvable* if for all $\epsilon < 1$, there exist $i, j \in \mathcal{A}$, for which

$$(10) \qquad \liminf_{n \to \infty} D_V(\mu_i^n[\epsilon], \mu_j^n[\epsilon]) > 0.$$

Note that by taking $\varepsilon = 0$ in (9) or (10) we obtain the original reconstruction condition (7). The same is true if $k = 0$ in (8).

Let $\lambda_2(M)$ denote the eigenvalue of $M$ which has the second largest absolute value [$\lambda_2(M)$ may be negative or nonreal]. In our main result we prove the following:

THEOREM 1.2. *Consider an ergodic Markov chain on the $B$-ary tree such that $B|\lambda_2(M)|^2 < 1$. Then we have the following:*

(i) *The reconstruction problem is not robust-solvable. Moreover, there exists $k^*$ such that for all $k > k^*$,*

$$(11) \qquad \max_{i,j} \lim_{n \to \infty} D_V(\mu_i^n[k], \mu_j^n[k]) = 0.$$

(ii) *For all nondegenerate $\nu$, the reconstruction problem is not $\nu$-robust-solvable. Moreover, for all nondegenerate $\nu$, there exists $\varepsilon^* < 1$ such that for all $\varepsilon > \varepsilon^*$,*

$$(12) \qquad \max_{i,j} \lim_{n \to \infty} D_V(\mu_i^n[\nu, \varepsilon], \mu_j^n[\nu, \varepsilon]) = 0.$$

(iii) *If all the entries of $M$ are nonzero, then the reconstruction problem is not erasure-robust-solvable. Moreover, there exists an $\epsilon^* < 1$ such that for all $\epsilon > \epsilon^*$,*

$$(13) \qquad \max_{i,j} \lim_{n \to \infty} D_V(\mu_i^n[\epsilon], \mu_j^n[\epsilon]) = 0.$$

It is easy to see that the total variation distances in (11), (12) and (13) are monotone decreasing in $k, \varepsilon$ and $\epsilon$ respectively.

The following proposition follows immediately from [12] or from the proofs in [22]. Together with Theorem 1.2, it shows that the threshold for robust reconstruction is given by $B|\lambda_2(M)|^2 = 1$.

PROPOSITION 1.3. *Consider an ergodic Markov chain on the $B$-ary tree where $B|\lambda_2(M)|^2 > 1$. Then we have the following:*

(i) *The reconstruction problem is robust-solvable.*
(ii) *For all nondegenerate $\nu$ the reconstruction problem is $\nu$-robust-solvable.*
(iii) *The reconstruction problem is erasure-robust-solvable.*

1.2. *Discussion.* The reconstruction problem was first studied in statistical physics [9, 28], where the problem was phrased in terms of extremality of the free measure for the Ising model on the $(B + 1)$-regular tree (Bethe lattice). It is not too hard to see (e.g., [3]) that the measure is nonextremal if and only if the reconstruction problem is solvable for the Markov chain on the $B$-ary tree with transition probabilities given by the binary symmetric Markov chain: $\begin{pmatrix} 1 - \delta & \delta \\ \delta & 1 - \delta \end{pmatrix}$. $\delta$ is related to the "inverse temperature" $\beta$ by $1 - 2\delta = \tanh(2\beta)$. The equivalence between nonextremality of the free measure of a random field and reconstruction solvability of an associated Markov chain on the same tree holds under mild nondegeneracy conditions (see, e.g., [18]).

In the past decade, the reconstruction problem reappeared in many applications: In communication networks (see [3] and the references there), in noisy computation (a model introduced by von Neumann in [29], see [8, 4]) and in phylogeny (molecular evolution, see [5, 26] for general background) [27]. Most recently, it is shown that the reconstruction problem is of crucial importance to basic questions in phylogeny [19, 20, 23]. In all of these applications the interest is to find when is it possible to reconstruct some information on the root state from states at the leaves of a finite tree. In many of the applications it is natural to consider robust-reconstruction as the observed data goes via additional "noise mechanism."

Solvability of the reconstruction problem is also closely related to the mixing rate of Glauber dynamics on the tree. See [1, 17], where it is shown that nonsolvability roughly corresponds to rapid mixing dynamics on the tree.

Determining if the reconstruction problem is solvable or not turns out to be very hard. Binary symmetric Markov chains is the only family for which the threshold for reconstruction solvability is known. Even here there is a generation gap between the proof of the lower-bound [9] and proofs of the upper bound [2] (see also [10] for a different proof, [3] for the result on general trees and [24] for the critical case on general trees). For binary symmetric Markov chains on the $B$-ary tree the threshold for the reconstruction problem is given by $B(1 - 2\delta)^2 = 1$, or, equivalently, $B\lambda_2(M)^2 = 1$. For all other families of Markov chain, including $q$-ary symmetric Markov chain for $q > 2$ and general $2 \times 2$ Markov chains, only bounds are known [15, 17, 22].

The threshold $B|\lambda_2(M)|^2 = 1$ is also the threshold for "census-solvability" [22], where different nodes of $L_n$ are indistinguishable (in other words, we only observe the "census" of level $n$). However, in general, it is not the threshold for reconstruction. Indeed, except for the binary symmetric channel, we know of no family of chains for which $B|\lambda_2(M)|^2 = 1$ is the threshold for reconstruction. Moreover, [18] shows that for asymmetric binary Markov chains (general stochastic 2 by 2 matrices) or symmetric Markov chains for $q > 2$ [where $M_{i,j} = (1-\delta)\mathbf{1}_{\{i=j\}} + \frac{\delta}{q-1}\mathbf{1}_{\{i \neq j\}}$ for $i, j \in \{1, \ldots, q\}$], the reconstruction problem is sometimes solvable even when $B|\lambda_2(M)|^2 < 1$. In [18] there is also a construction of $M$ with $\lambda_2(M) = 0$ for which the reconstruction problem is solvable for large $B$.

Why is determining the threshold for reconstruction hard? From the technical point of view a Markov chain on the tree corresponds to a recursion in some random variables ([2, 17, 24]). A natural way to analyze these recursions is to use a perturbative argument around the stationary distribution of the chain. The main problem is that the random variables we start with are atoms—far from the stationary distribution—and that, in general, the recursions lack any convexity. For "robust-reconstruction" the problem is easier—as the recursions begin close to the stationary distribution.

Our proof is based on a new measure of *discrepancy* for a vector of distributions which is a weighted variant of the $\chi^2$ distance. We show that an application of the chain $M$ contracts the discrepancy by a $|\lambda_2(M)|^2$ factor, and that if the discrepancy is smaller than $\delta$, then tensoring $B$ copies of the distributions increases the discrepancy by a factor of at most $B(1 + \varepsilon(\delta))$, where $\varepsilon(\delta) \to 0$, as $\delta \to 0$.

It is interesting to compare our results with the results of [25]. In [25] Pemantle and Steif study robust phase transition on trees. For a Gibbs measure on a tree we say that a robust phase transition occurs if the boundary conditions on a cutset have a nonvanishing effect on the root even when the interactions along the cutsets are made arbitrarily small but fixed (see [25] for exact definition). The main results of [25] give the exact threshold for robust phase transitions for general (bounded degree) trees for Potts and Heisenberg models in terms of the underlying model and the branching number (see [14]) of the tree.

Both in our result and in the results of [25], it is easier to analyze the "robust" problem than it is the original problem for similar reasons. In both cases the "nonrobust" problem is hard to control without some convexity assumption, while the solution of "robust" problem allows the use of "local" arguments.

Moreover, like robust phase transition, robust reconstruction is a geometric property, that is, for general bounded degree $T$, the threshold for robust-reconstruction depends only on $\text{br}(T)$ and $|\lambda_2(M)|$. Indeed, the proof of Theorem 3.3 combines the analysis of the new discrepancy measure introduced here, with some of the techniques developed in [25] for controlling recursions on general trees.

A natural open problem is to determine the behavior of robust-solvability in the critical case, where $B|\lambda_2(M)|^2 = 1$. Our techniques shed no light on this problem. It is also interesting to try and remove the restriction that the entries of $M$ are positive for (13); see also Remark 2.10. Finally, in the proof presented for Theorem 3.3, for fixed $M$ and $\nu$, the bounds on $\varepsilon$ and $k$ are becoming weaker as $\text{br}(T)|\lambda_2(M)|^2$ approaches 1 (i.e., $\varepsilon \to 1$ and $k \to \infty$). It is natural to ask if for given $M$, $\nu$ and $K$, there exist $\varepsilon$ and $k$ for which the result holds uniformly for all infinite trees $T$ with $\text{br}(T)|\lambda_2(M)|^2 < 1$.

**2. Proof.** Recall that we denote by $M$ the transition matrix. In this section we will often multiply $M$ from the right by a vector of functions, from the left by a vector of measures—in which case the resulting vector would also be a vector of functions/vector of measures.

Let $\mathbf{1} = (1, \ldots, 1)^t$, then clearly $M\mathbf{1} = \mathbf{1}$. Let $v_i$ be the stationary probability of state $i$, and $v = (v_1, \ldots, v_q)$ the stationary distribution, so that $vM = v$. In the remainder of this section we will use $a, b, c, \ldots$ for column vectors, and $u, v, w, \ldots$ for row vectors. $v$ will always denote the stationary distribution.

Note that if $b$ is a column vector such that $vb = 0$, then $vMb = vb = 0$. In other words, the linear space $v^\perp = \{b \in \mathbb{R}^q : vb = 0\}$ is invariant under $M$.

LEMMA 2.1. *Let $b = (b_1, \ldots, b_q)^t$ be a vector such that $b_i > 0$ for all $i$. Then*

$$\sum_{i=1}^{q} \frac{v_i}{(Mb)_i} \leq \sum_{i=1}^{q} \frac{v_i}{b_i}.$$

PROOF. By Jensen,

$$\frac{1}{(Mb)_i} = \frac{1}{\sum_j m_{i,j} b_j} \leq \sum_j m_{i,j} \frac{1}{b_j}.$$

Hence,

$$\sum_i \frac{v_i}{(Mb)_i} \leq \sum_i \sum_j v_i m_{i,j} \frac{1}{b_j} = \sum_j \left( \sum_i v_i m_{i,j} \right) \frac{1}{b_j} = \sum_j \frac{v_j}{b_j},$$

as needed. □

By looking at the Jordan form of $M$ it is easy to see the following:

Lemma 2.2.  *Given $\varepsilon > 0$, there exists an Euclidean norm $\|\cdot\|$ on $v^{\perp}$ such that $\|Mb\| \leq (|\lambda_2(M)| + \varepsilon)\|b\|$ for all $b \in v^{\perp}$.*

Let $Q$ be the projection onto $v^{\perp}$ defined by $Qb = b - (vb)\mathbf{1}$ [note that $vQb = vb - (vb)(v\mathbf{1}) = 0$ for all $b$].

Definition 2.3.  Let $\|\cdot\|$ be a Euclidean norm on $v^{\perp}$. Let $\nu = (\nu_1, \ldots, \nu_q)$ be a vector of distributions on a common space. Let $f = (f_1, \ldots, f_q)$ be the vector of density functions with respect to a $\sigma$-finite measure $\mu$, such that $\nu_i \ll \mu$ for every $i$. In other words, $d\nu_i = f_i\, d\mu$ for all $i$. We then define the *discrepancy* of the vector by

$$D_{\|\cdot\|}^{\mu}(f) = \int \|Qf\|^2 \sum_{i=1}^{q} \frac{v_i}{f_i}\, d\mu.$$

We also write $D(\nu)$ and $D(f)$ for the discrepancy, without explicitly indicating the norm and reference measure.

Note the similarity between the discrepancy and the $\chi^2$-distance. The $\chi^2$-distance is known to be well behaved with respect to $\ell_2$ norms. Note that if $f_1 = \cdots = f_q$, then $Qf = fQ\mathbf{1} = 0$. Thus, $Q$ projects into the orthogonal complement of the space where the discrepancy should be 0.

Lemma 2.4.  *$D(\nu)$ is independent of the reference measure $\mu$; that is, if $f = (f_1, \ldots, f_q)$ and $g = (g_1, \ldots, g_q)$ are such that $d\nu_i = f_i\, d\mu = g_i\, d\tilde{\mu}$ for all $i$, then $D_{\|\cdot\|}^{\mu}(f) = D_{\|\cdot\|}^{\tilde{\mu}}(g)$.*

Proof.  Assume that $\tilde{\mu} \ll \mu$. The general case then follows by considering the three reference measures $\mu$, $\mu + \tilde{\mu}$, $\tilde{\mu}$.

Since $Q$ is linear and $\|\cdot\|$ is Euclidean, we may write $\|Qb\|^2 = \sum_{i,j=1}^{q} t_{i,j}b_i b_j$ for some $t_{i,j} \in \mathbb{R}$. Now

$$D^{\mu}(f) = \sum_{r,i,j=1}^{q} v_r t_{i,j} \int \frac{f_i f_j}{f_r}\, d\mu = \sum_{r,i,j=1}^{q} v_r t_{i,j} \int \frac{g_i g_j (d\tilde{\mu}/d\mu)^2}{g_r\, d\tilde{\mu}/d\mu}\, d\mu$$

$$= \sum_{r,i,j=1}^{q} v_r t_{i,j} \int \frac{g_i g_j}{g_r}\, d\tilde{\mu} = D^{\tilde{\mu}}(g). \qquad \square$$

Lemma 2.5.  *Let the norm $\|\cdot\|$ on $v^{\perp}$ satisfy $\|Mb\| \leq \alpha\|b\|$ for all $b \in v^{\perp}$ and some constant $\alpha$. Then $D(Mf) \leq \alpha^2 D(f)$ for all $f$, as in Definition 2.3.*

Proof.  For all $b$,

$$MQb = M(b - (vb)\mathbf{1}) = Mb - (vb)M\mathbf{1}$$
$$= Mb - (vb)\mathbf{1} = Mb - (vMb)\mathbf{1} = QMb.$$

Therefore, we have pointwise that

$$(14) \qquad \|QMf\|^2 = \|MQf\|^2 \le \alpha^2 \|Qf\|^2.$$

Now

$$D(Mf) = \int \|QMf\|^2 \sum_{i=1}^q \frac{v_i}{(Mf)_i} \, d\mu$$

$$(15) \qquad \le \alpha^2 \int \|Qf\|^2 \sum_{i=1}^q \frac{v_i}{(Mf)_i} \, d\mu$$

$$(16) \qquad \le \alpha^2 \int \|Qf\|^2 \sum_{i=1}^q \frac{v_i}{f_i} \, d\mu = \alpha^2 D(f),$$

where (15) follows from (14), and (16) follows by Lemma 2.1. $\square$

LEMMA 2.6. *For every Euclidean norm $\|\cdot\|$ on $v^\perp$, there exists a constant $C = C(\|\cdot\|)$, such that*

$$(17) \qquad \left| \int \frac{f_i f_j}{f_k} \, d\mu - 1 \right| \le C D(f),$$

*for all $i, j, k$, where $f$ and $\mu$ are as in Definition 2.3.*

PROOF. By Cauchy–Schwarz,

$$\left| \int \frac{f_i f_j}{f_k} \, d\mu - 1 \right| = \left| \int \frac{(f_i - f_k)(f_j - f_k)}{f_k} \, d\mu \right|$$

$$\le \sqrt{\int \frac{|f_i - f_k|^2}{f_k} \, d\mu} \sqrt{\int \frac{|f_j - f_k|^2}{f_k} \, d\mu}.$$

Therefore, in order to prove the lemma, it suffices to prove that there exists a constant $C$ such that for all $i, j, k$, it holds that

$$(18) \qquad \int \frac{|f_i - f_j|^2}{f_k} \, d\mu \le C D(f).$$

Note that

$$(Qb)_i - (Qb)_j = (b - (vb)\mathbf{1})_i - (b - (vb)\mathbf{1})_j = b_i - b_j.$$

Therefore, for all $b$, it holds that, for some constant $C_{ij}$,

$$|b_i - b_j| = |(Qb)_i - (Qb)_j| \le C_{ij} \|Qb\|.$$

Hence,

$$\int \frac{|f_i - f_j|^2}{f_k} \, d\mu \le C_{ij}^2 \int \frac{\|Qf\|^2}{f_k} \le \frac{C_{ij}^2}{v_k} D(f).$$

Now (18) follows by taking $C = \sup_{i,j,k} \frac{C_{ij}^2}{v_k}$. $\square$

Given a $\sigma$-finite measure $\mu$ on a space $X$, we denote by $\mu^{\otimes B}$ the product measure on $X^B$ with marginals $\mu$. Similarly, if $f_i$ is a density of $\nu_i$ with respect to $\mu$, write $f_i^{\otimes B}$ for the density of $\nu_i^{\otimes B}$ with respect to $\mu^{\otimes B}$. Finally, for $f = (f_1, \ldots, f_q)$, write $f^{\otimes B} = (f_1^{\otimes B}, \ldots, f_q^{\otimes B})$, and for $\nu = (\nu_1, \ldots, \nu_q)$, write $\nu^{\otimes B} = (\nu_1^{\otimes B}, \ldots, \nu_q^{\otimes B})$. We similarly use $\bigotimes_{r=1}^B \nu^r$ for the componentwise product of several vectors $\nu^r$ of measures.

LEMMA 2.7. *Let $\|\cdot\|$ be an Euclidean norm on $v^\perp$, $B \geq 1$ an integer and $\varepsilon > 0$. Then there exists a $\delta > 0$ such that if $\nu^1 = (\nu_1^1, \ldots, \nu_q^1)$, ..., $\nu^B = (\nu_1^B, \ldots, \nu_q^B)$, satisfy $D(\nu^i) \leq \delta$ for $1 \leq i \leq B$, then*

$$D\left(\bigotimes_{r=1}^B \nu^r\right) \leq (1+\varepsilon)(D(\nu^1) + \cdots + D(\nu^B)).$$

*In particular, given $\varepsilon > 0$, there exists a $\delta > 0$ such that if $\nu = (\nu_1, \ldots, \nu_q)$ satisfies $D(\nu) \leq \delta$, then $D(\nu^{\otimes B}) \leq (1+\varepsilon)BD(\nu)$.*

PROOF. The second part of the lemma immediately follows from the first part. Choose a reference measure $\mu$ with $\nu_i^r \ll \mu$ for every $i$ and $r$, and let $f_i^r$ be the density $d\nu_i^r/d\mu$.

As in Lemma 2.4, we may write $\|Qb\|^2 = \sum_{i,j} t_{i,j} b_i b_j$. Moreover, since $\|Q\mathbf{1}\|^2 = \|0\|^2 = 0$, it follows that $\sum_{i,j} t_{i,j} = 0$. Hence,

$$
\begin{aligned}
D(f) &= \int \|Qf\|^2 \sum_{k=1}^q \frac{v_k}{f_k} \, d\mu \\
&= \sum_{i,j,k} v_k t_{i,j} \int \frac{f_i f_j}{f_k} \, d\mu \\
&= \sum_{i,j,k} v_k t_{i,j} \left(\int \frac{f_i f_j}{f_k} \, d\mu - 1\right).
\end{aligned}
$$
(19)

Substituting $\bigotimes_{r=1}^B f^r$ in (19), we obtain, using the reference measure $\mu^{\otimes B}$,

$$
\begin{aligned}
D\left(\bigotimes_{r=1}^B f^r\right) &= \sum_{i,j,k} v_k t_{i,j} \left(\int \frac{\bigotimes_{r=1}^B f_i^r \bigotimes_{r=1}^B f_j^r}{\bigotimes_{r=1}^B f_k^r} \, d\mu^{\otimes B} - 1\right) \\
&= \sum_{i,j,k} v_k t_{i,j} \left(\prod_{r=1}^B \left(\int \frac{f_i^r f_j^r}{f_k^r} \, d\mu\right) - 1\right).
\end{aligned}
$$
(20)

Let $C$ be chosen to satisfy (17) in Lemma 2.6, and $\tilde{C} = C \sum_{i,j} |t_{i,j}|$. Let $\delta$ be chosen such that for all $(x_1, \ldots, x_B) \in [1 - C\delta, 1 + C\delta]^B$, it holds that

$$(21) \qquad \left| \prod_{r=1}^{B} x_r - 1 - \sum_{r=1}^{B}(x_r - 1) \right| \leq \frac{\varepsilon \sum_{r=1}^{B} |x_r - 1|}{\tilde{C}}.$$

By Lemma 2.6, it follows that if $D(f) \leq \delta$, then for all $i, j, k$,

$$\left| \int \frac{f_i f_j}{f_k} \, d\mu - 1 \right| \leq C D(f) \leq C\delta.$$

Therefore, it follows from (20) and (21) that

$$
\begin{aligned}
(22) \qquad D\left( \bigotimes_{r=1}^{B} f^r \right) &\leq \sum_{i,j,k,r} v_k t_{i,j} \left( \int \frac{f_i^r f_j^r}{f_k^r} \, d\mu - 1 \right) + \frac{\varepsilon}{\tilde{C}} \sum_{i,j,k,r} v_k |t_{i,j}| \left| \int \frac{f_i^r f_j^r}{f_k^r} \, d\mu - 1 \right| \\
&= \sum_r D(f^r) + \frac{\varepsilon}{\tilde{C}} \sum_{i,j,k,r} v_k |t_{i,j}| \left| \int \frac{f_i^r f_j^r}{f_k^r} \, d\mu - 1 \right| \\
&\leq \sum_r D(f^r) + \frac{\varepsilon}{\tilde{C}} \sum_{i,j,k,r} v_k |t_{i,j}| C D(f^r) \\
&= (1 + \varepsilon) \sum_r D(f^r),
\end{aligned}
$$

where inequality (22) follows from Lemma 2.6. $\quad \square$

LEMMA 2.8. *Given a Euclidean norm $\|\cdot\|$ on $v^\perp$, there exists a constant $C(\|\cdot\|) < \infty$ such that for any vector $\nu = (\nu_i)_{i=1}^q = (f_i \, d\mu)_{i=1}^q$ of distributions we have*

$$(23) \qquad \sup_{i,j} d_{\mathrm{TV}}(\nu_i, \nu_j) = \sup_{i,j} \int |f_i - f_j| \, d\mu \leq C \sqrt{D(f)},$$

*where $f = (f_1, \ldots, f_q)$.*

PROOF. By Cauchy–Schwarz,

$$\int |f_i - f_j| \, d\mu \leq \sqrt{\int \frac{|f_i - f_j|^2}{f_i} \, d\mu} \sqrt{\int f_i \, d\mu} = \sqrt{\int \frac{|f_i - f_j|^2}{f_i} \, d\mu},$$

and (23) follows from Lemma 2.6. $\quad \square$

LEMMA 2.9. *(i) Let $\|\cdot\|$ be a Euclidean norm on $v^\perp$. Let $\mu$ be a probability distribution on $1, \ldots, q$ such that $\mu(i) > 0$ for all $i$. Then for all $\delta > 0$, there exists an $\varepsilon^* = \varepsilon^*(\delta) < 1$ such that for any vector $\nu' = (\nu'_1, \ldots, \nu'_q)$ of probability distributions on $1, \ldots, q$ and for all $\varepsilon > \varepsilon^*$, if*

$$\nu = (1 - \varepsilon)(\nu'_1, \ldots, \nu'_q) + \varepsilon(\mu, \ldots, \mu),$$

*then $D(\nu) \leq \delta$.*

(ii) *Let* $\|\cdot\|$ *be a Euclidean norm on* $v^{\perp}$. *Let* $\nu_{\ell}^{(r)}$ *denote the $\ell$th row of* $M^r$, *that is, the distribution of the chain given by $M^r$, after $r$ steps starting at $\ell$, and $\nu^{(r)} = (\nu_1^{(r)}, \ldots, \nu_q^{(r)})$. Then for all $\delta > 0$, there exists an $r^* = r^*(\delta)$ such that if $r \geq r^*$, then $D(\nu^{(r)}) \leq \delta$.*

(iii) *Let* $\|\cdot\|$ *be an Euclidean norm on* $v^{\perp}$. *Let* $\nu = (\nu_1, \ldots, \nu_q)$ *be a vector of probability measures such that for all $1 \leq i, j \leq q$, it holds that $\nu_i(j) > 0$. Let $\nu^{\epsilon} = (\nu_1^{\epsilon}, \ldots, \nu_q^{\epsilon})$ be a collection of probability measures on $\{1, \ldots, q+1\}$ such that $\nu_i^{\epsilon} = (1-\epsilon)\nu_i + \epsilon\nu'$, where $\nu'$ is the delta measure on $q+1$. Then for all $\delta > 0$, there exists an $\epsilon^* = \epsilon^*(\delta) < 1$ such that if $\epsilon > \epsilon^*$, then $D(\nu^{\epsilon}) \leq \delta$.*

Proof. For the first part of the lemma we use the representation of $D(\nu)$ as in (19) with respect to the measure $\mu$. Let $m = \min\{\mu(1), \ldots, \mu(q)\}$, and observe that if $d\nu_i = f_i \, d\mu$, then

$$f_i = \frac{d\nu_i}{d\mu} = \varepsilon + (1-\varepsilon)\frac{d\nu'_i}{d\mu},$$

so

$$\varepsilon \leq f_i \leq \varepsilon + (1-\varepsilon)/m,$$

and for all $i, j, k$,

$$\frac{f_i f_j}{f_k} \leq \frac{(\varepsilon + (1-\varepsilon)/m)^2}{\varepsilon}.$$

Hence, by (19),

$$(24) \qquad D(\nu) \leq \left(\frac{(\varepsilon + (1-\varepsilon)/m)^2}{\varepsilon} - 1\right)\sum_{i,j}|t_{i,j}|,$$

and the right-hand side of (24) converges to 0 as $\varepsilon \to 1$.

The second part of the lemma follows from the first one, as the ergodicity of $M$ implies that for all $i$, $\nu_i^{(m)}$ converges to the stationary distribution of the chain as $m \to \infty$.

The third part of the lemma is proven similarly to the first part. Let $m = \min_{i,j \in \mathcal{A}} \nu_i(j)$ and

$$\mu^{\varepsilon} = \varepsilon\nu' + \frac{(1-\varepsilon)}{q}\sum_{i=1}^{q}\nu_i.$$

Note that if $d\nu_i^{\varepsilon} = f_i \, d\mu^{\varepsilon}$, then $m \leq f_i(\ell) \leq q$ for $1 \leq \ell \leq q$ and all $i$, and $f_i(q+1) = 1$ for all $i$. It follows that for all $i, j, k$ and $1 \leq \ell \leq q$,

$$\frac{m^2}{q} \leq \frac{f_i(\ell)f_j(\ell)}{f_k(\ell)} \leq \frac{q^2}{m} \quad \text{and} \quad \frac{f_i(q+1)f_j(q+1)}{f_k(q+1)} = 1.$$

Moreover, $\mu^\varepsilon(q+1) = \varepsilon$ and $\mu^\varepsilon(\{1,\ldots,q\}) = 1 - \varepsilon$, so for all $i, j$ and $k$,

$$\int \frac{f_i f_j}{f_k} \, d\mu - 1 = \int \left( \frac{f_i f_j}{f_k} - 1 \right) d\mu \leq \left( \frac{q^2}{m} - 1 \right)(1 - \varepsilon).$$

Hence, by (19),

$$(25) \qquad D(\nu^\varepsilon) = \sum_{i,j,k} v_k t_{i,j} \left( \int \frac{f_i f_j}{f_k} \, d\mu - 1 \right) \leq \left( \frac{q^2}{m} - 1 \right)(1 - \varepsilon) \sum_{i,j} |t_{i,j}|$$

and the right-hand side of (25) converges to 0 as $\varepsilon \to 1$. $\quad\square$

PROOF OF THEOREM 1.2. The basic idea of the proof is that if $\mu[N]^n$ is the vector of probability measures $(\mu[N]_1^n, \ldots, \mu[N]_q^n)$ defined in (5), then we may write $\mu[N]^{n+1}$ in terms of $\mu[N]^n$ using the operator $M$ and tensoring. This will allow us to bound discrepancies recursively. Let $\rho_1, \ldots, \rho_B$ be the $B$ children of $\rho$ in the $B$-ary tree. Write $E_{n+1}(s)$ for the edges in $E_{n+1}$ that are on the subtree rooted in $\rho_s$ (formally, these are the edges of $E_{n+1}$ that are connected to $\rho$ only by paths going via $\rho_s$). Define $L_{n+1}(s)$ similarly. Finally, for a configuration $\sigma$ of the vertices at the first $n$ level of the tree, let $\sigma_1, \ldots, \sigma_B$ denote the configurations restricted to the subtrees rooted at $\rho_1, \ldots, \rho_B$. Then by (4) and (5),

$$
\begin{aligned}
\mu[N]_\ell^{n+1}(\tau) &= \sum_\sigma \mathbf{1}_{\{\sigma(\rho)=\ell\}} \prod_{(x,y) \in E_{n+1}} M_{\sigma(x),\sigma(y)} \prod_{y \in L_{n+1}} N_{\sigma(y),\tau(y)} \\
&= \prod_{s=1}^{B} \left( \sum_{\ell'=1}^{q} M_{\ell,\ell'} \sum_{\sigma_s} \mathbf{1}_{\{\sigma_s(\rho_s)=\ell'\}} \prod_{(x,y) \in E_{n+1}(s)} M_{\sigma_s(x),\sigma_s(y)} \right. \\
&\qquad\qquad\qquad\qquad\qquad \left. \times \prod_{y \in L_{n+1}(s)} N_{\sigma_s(y),\tau_s(y)} \right) \\
&= \prod_{s=1}^{B} \left( \sum_{\ell'=1}^{q} M_{\ell,\ell'} \mu[N]_{\ell'}^n(\tau | L_{n+1}(s)) \right).
\end{aligned}
$$
(26)

Note that the expression in the parenthesis in (26) is given by $(M(\mu[N]^n))_\ell(\tau | L_{n+1}(s))$, the $\ell$th coordinate of the vector $M(\mu[N]^n)$. It is now easy to see that

$$(27) \qquad \mu[N]^{n+1} = (M(\mu[N]^n))^{\otimes B}.$$

We use (27) in order to bound discrepancies recursively.

The assumption $B|\lambda_2(M)|^2 < 1$ implies by Lemma 2.2 that there exists an $\varepsilon > 0$ and a norm $\|\cdot\|$ on $v^\perp$ such that for all $b \in v^\perp$, it holds that $\|Mb\| \leq \alpha\|b\|$, where $B(1+\varepsilon)\alpha^2 \leq 1 - \varepsilon$. Let $\delta$ be chosen as to satisfy Lemma 2.7, so that if $D(f) \leq \delta$, then $D(f^{\otimes B}) \leq B(1+\varepsilon)D(f)$.

By Lemma 2.9 it follows that there exists a $k^*$ such that $\mu^0 = (\mu_1^0[k], \ldots, \mu_q^0[k])$ satisfies $D(\mu^0) \leq \delta$ for all $k \geq k^*$. Write $\mu^n$ for $(\mu_1^n[k], \ldots, \mu_q^n[k])$. Thus, (27) implies that $\mu^{n+1} = (M\mu^n)^{\otimes B}$. It now follows by Lemmas 2.7 and 2.5 that

$$D(\mu^{n+1}) \leq B(1+\varepsilon)\alpha^2 D(\mu^n) \leq (1-\varepsilon)D(\mu^n).$$

Hence, $\lim_{n\to\infty} D(\mu^n) = 0$. We therefore conclude from Lemma 2.8 that

$$\lim_{n\to\infty} \max_{i,j} d_{\mathrm{TV}}(\mu_i^n[k], \mu_j^n[k]) = 0,$$

and (11) follows.

In order to prove (12), let $\nu$ be a nondegenerate measure and note that by Lemma 2.9 it follows that there exist an $\varepsilon^* < 1$ such that $\mu^0 = (\mu_1^0[\nu, \varepsilon], \ldots, \mu_q^0[\nu, \varepsilon])$ satisfies $D(\mu^0) \leq \delta$ for all $\varepsilon > \varepsilon^*$. Now (12) follows similarly to (11).

The proof of (13) is similar. We look at $\mu^1 = (\mu_1^1[\epsilon], \ldots, \mu_q^1[\epsilon])$. Note that (27) implies that $\mu^1 = (M\mu^0)^{\otimes B}$. Let $\nu = M\mu^0$. Note that the vector $\nu$ satisfies for all $i$ that $\nu_i(q+1) = \epsilon$ and $\nu_i(j) = (1-\epsilon)M_{i,j}$ otherwise.

Since all the entries of $M$ are positive, it follows from Lemma 2.9 that for every $\delta' > 0$ there exists an $\epsilon^* < 1$ such that if $\epsilon > \epsilon^*$, then $D(\nu) \leq \delta'$. We may now apply Lemma 2.7 and choose $\delta' > 0$ in such a way that

$$D(\mu^1) = D(\nu^{\otimes B}) \leq \delta.$$

The rest of the proof is identical.  □

REMARK 2.10. It is an interesting goal to extend (13) to general ergodic chains (where some of the entries of $M$ may be zero). Above we proved this for the case where all the entries of $M$ are positive.

The proof of Lemma 2.9 can easily be extended to the case when there exists an $n$ such that the measures $\mu_1^n, \ldots, \mu_q^n$ all have the same support [in such a case one can prove that there exists a value of $\epsilon < 1$ such that $D(\mu^n) \leq \delta$ by showing that for $\epsilon$ sufficiently large, the measures $\mu_i^n$ have most of their mass on the atom $(q+1, \ldots, q+1)$ and bounded relative densities elsewhere]. However, we do not know any simple characterizations of the matrices $M$ for which this holds (it evidently depends only on the set of zero entries of $M$); nor we do believe that this property is necessary for nonerasure-robust-solvability.

## 3. General trees.

Our results readily extend to general infinite bounded degree trees, where $B$ is replaced by $\mathrm{br}(T)$, the branching number of the tree. In [6], Furstenberg introduced the notion of the Hausdorff dimension of a tree. Later, Lyons [13, 14] showed how many of the probabilistic properties of the tree are determined by this number which he named the branching number.

For our purposes it is best to define the branching number via cutsets. A *cutset* $S$ for a tree $T$ rooted at $\rho$ is a finite set of vertices separating $\rho$ from $\infty$. In other words, a finite set $S$ is a cutset, if every infinite self avoiding path from $\rho$ intersects $S$. An *antichain* is a cutset that does not have any proper subset which is also a cutset.

We follow the notation of [25] and for a cutset $S$, write $\text{Ins}(S)$ for the inside of $S$ (the finite component of $T \setminus S$, containing the root $\rho$), $\text{Ins}_E(S)$ for edges inside $S$ [those edges $(x,y)$ having $x \in \text{Ins}(S)$] and $\text{Out}(S)$ for the outside of $S$ [$\text{Out}(S) = T \setminus (S \cup \text{Ins}(S))$].

DEFINITION 3.1. The branching number of $T$ is defined as

$$\text{br}(T) = \inf \left\{ \lambda > 0 \colon \inf_{\text{cutsets } S} \sum_{x \in S} \lambda^{-|x|} = 0 \right\}.$$

Note that

$$\text{br}(T) = \inf \left\{ \lambda > 0 \colon \inf_{\text{antichains } S} \sum_{x \in S} \lambda^{-|x|} = 0 \right\}.$$

By Min-Cut-Max-Flow, $\text{br}(T)$ is also the supremum of the real numbers $\lambda > 0$, such that $T$ admits a positive flow from the root to infinity, where on every edge $e$ of $T$, the flow is bounded by $\lambda^{-|e|}$. It is shown in [14] that $\text{br}(T)^{-1}$ is the critical probability for Bernoulli percolation on $T$. See [14] and [3] for equivalent definitions of $\text{br}(T)$ in terms of percolation, cutset sums and electrical conductance. We note that $\text{br}(T_B) = B$ for the $B$-ary tree $T_B$.

As in Section 1, the Markov chain on $T$ is described by an $|\mathcal{A}| \times |\mathcal{A}|$ stochastic matrix $M$ and the perturbations by an $|\mathcal{A}| \times |\mathcal{B}|$ stochastic matrix $N$. For $B$-ary trees, we observed the process on the special antichains $L_n$; for general trees, it seems more natural to consider arbitrary antichains. The distribution $\mu[N]_\ell^S$ of the observed (perturbed) states on an antichain $S$ in $T$ is given by, extending (4) and (5),

$$(28) \quad \mu[N]_\ell^S(\tau) = \sum_\sigma \mathbf{1}_{\{\sigma(\rho) = \ell\}} \prod_{(x,y) \in \text{Ins}_E(S)} M_{\sigma(x),\sigma(y)} \times \prod_{y \in S} N_{\sigma(y),\tau(y)}.$$

We proceed by defining $\mu_\ell^S$, the measure $\mu_\ell^S[k]$ for $k \geq 0$, the measure $\mu_\ell^S[\nu,\varepsilon]$ for $\varepsilon > 0$ and nondegenerate distribution $\nu$ on $\mathcal{A}$ and $\mu_\ell^S[\epsilon]$. This is done in exactly the same way as in the case of the $B$-ary tree, by choosing appropriate $N$'s in (28).

We say that the reconstruction problem is solvable if there exists $i,j \in \mathcal{A}$, for which

$$\inf_{S \text{ antichain}} D_V(\mu_i^S, \mu_j^S) > 0,$$

where $\mu_\ell^S$ denotes the conditional distribution on $\sigma_S$ given that $\sigma(\rho) = \ell$. We similarly define the notions of robust-solvable, $\nu$-robust-solvable and erasure-robust-solvable.

REMARK 3.2. The definitions of solvability for general trees and $B$-ary tree are *not* compatible. If $T$ is the $B$-ary tree, then solvability by Definition 1.1 involves only cutsets $S = L_n$ and is therefore a weaker condition than solvability defined here, which involves all antichains (same for robust-solvable etc.). However, we will obtain the same threshold for robust-reconstruction under both definitions.

The proof of our main result extends to show the next theorem.

THEOREM 3.3. *Consider an ergodic Markov chain on a rooted tree $T$ such that $\mathrm{br}(T)|\lambda_2(M)|^2 < 1$. Then we have the following:*

  (i) *The reconstruction problem is not robust-solvable.*
  (ii) *For all nondegenerate $\nu$, the reconstruction problem is not $\nu$-robust-solvable.*
  (iii) *If all the entries of $M$ are nonzero, then the reconstruction problem is not erasure-robust-solvable.*

This proves that the threshold for robust reconstruction is given by $\mathrm{br}(T) \times |\lambda_2(M)|^2 = 1$ as the proof of Theorem 1.4 in [22] immediately generalizes to show the following proposition:

PROPOSITION 3.4. *Consider an ergodic Markov chain on a tree $T$ such that $\mathrm{br}(T)|\lambda_2(M)|^2 > 1$. Then the reconstruction problem is robust-solvable, for all nondegenerate $\nu$ the reconstruction problem is $\nu$-robust-solvable and the reconstruction problem is erasure-robust-solvable.*

We now turn to the proof of Theorem 3.3 which generalizes the proof of Theorem 1.2. For a vertex $x$ of the rooted tree $T$, we write $T(x)$ for the subtree rooted at $x$, that is, the subtree consisting of $x$ and all of its descendents. We will use the following lemma from Pemantle and Steif [25].

LEMMA 3.5 ([25], Lemma 3.3). *Assume that $\mathrm{br}(T) < g$. Then for all $\varepsilon > 0$ there exists an anitchain $S$ such that*

$$\tag{29} \sum_{x \in S} \left( \frac{1}{g} \right)^{|x|} \le \varepsilon,$$

*and for all $y \in S \cup \mathrm{Ins}(S)$,*

$$\tag{30} \sum_{x \in S \cap T(v)} \left( \frac{1}{g} \right)^{|x|-|y|} \le 1.$$

PROOF OF THEOREM 3.3. We will show that under the conditions of Theorem 3.3 the reconstruction problem is not robust-solvable.

Let $S$ be an antichain and $y \in S \cup \text{Ins}(S)$. Consider the Markov chain on the subtree $T(y)$, starting with state $\ell$ at $y$, and let $\mu[N]_\ell^{y,S}$ be the distribution of the observed states on $T(y) \cap S$. Thus,

$$(31) \quad \mu[N]_\ell^{y,S}(\tau) = \sum_\sigma \mathbf{1}_{\{\sigma(y)=\ell\}} \prod_{\substack{(x,z)\in\text{Ins}_E(S) \\ x\in T(y), z\in T(y)}} M_{\sigma(x),\sigma(z)} \times \prod_{z\in T(y)\cap S} N_{\sigma(z),\tau(z)}.$$

We write $\mu[N]^{y,S}$ for $(\mu[N]_1^{y,S},\dots,\mu[N]_q^{y,S})$. Note that $\mu[N]^S = \mu[N]^{\rho,S}$ for any antichain $S$.

The proof of Theorem 1.2 easily extends to show that if $y \in \text{Ins}(S)$ and $z_1,\dots,z_B$ are the children of $y$, then

$$(32) \quad \mu^{y,S}[N] = \bigotimes_{r=1}^B (M\mu^{z_r,S}[N]).$$

We will prove the theorem by recursively analyzing discrepancies via (32). We will prove the result for robust-solvability and indicate the modifications needed for other cases at the end of the proof.

Below, we will write $\mu^{y,S}[k]$ for the measure $\mu^{y,S}[N]$, where $N = M^k$. Note that if $S$ is an antichain and $v \in S$, then the measure $\mu^{y,S}$ is a measure on a single node. We may therefore apply Lemma 2.9 and conclude that for all $\delta > 0$, there exists a $k^*$ such that for $k \geq k^*$, for all antichains $S$ and $y \in S$, it holds that

$$(33) \quad D(\mu^{y,S}[k]) \leq \delta.$$

Since $\text{br}(T)|\lambda_2(M)|^2 < 1$, there exist, by Lemma 2.2, an $\varepsilon > 0$ and a norm $\|\cdot\|$ on $v^\perp$ such that for all $b \in v^\perp$, it holds that $\|Mb\| \leq \alpha\|b\|$, where $(1+\varepsilon)\,\text{br}(T)\alpha^2 \leq 1-\varepsilon$. Recall that there is a uniform bound $K$ on the number of children of vertices of $T$. Let $\delta$ be chosen as to satisfy Lemma 2.7 for every $B \leq K$, so that if $D(f^r) \leq \delta$ for $r = 1,\dots,B$, and $B \leq K$, then $D(\bigotimes_{r=1}^B f^r) \leq (1+\varepsilon)\sum_{r=1}^B D(f^r)$.

Lemma 3.5 implies that there exists a sequence of antichains $S_n$ such that

$$(34) \quad \lim_{n\to\infty} \sum_{x\in S_n} [(1+\varepsilon)\alpha^2]^{|x|} = 0,$$

and that for all $n$ and $y \in S_n \cup \text{Ins}(S_n)$,

$$(35) \quad \sum_{x\in S_n\cap T(y)} [(1+\varepsilon)\alpha^2]^{|x|-|y|} \leq 1.$$

We will now show by induction (on $s - |y|$, where $s = \max_{x \in S} |x|$), that for all antichains $S = S_n$ and all $y \in S \cup \mathrm{Ins}(S)$, for $k \geq k^*$,

$$(36) \qquad D(\mu^{y,S}[k]) \leq \delta \sum_{x \in S \cap T(y)} [(1 + \varepsilon)\alpha^2]^{|x| - |y|}.$$

The case where $y \in S$ follows from (33). This also proves the base of the induction. For the induction step, it therefore suffices to consider $y \in \mathrm{Ins}(S)$ such that the children of $v$ denoted $z_1, \ldots, z_B$ satisfy the induction hypothesis. By Lemma 2.5 and the induction hypothesis, for all $r$,

$$(37) \qquad D(M\mu^{z_r,S}[k]) \leq \alpha^2 D(\mu^{z_r,S}[k]) \leq \delta\alpha^2 \sum_{x \in S \cap T(z_r)} [(1 + \varepsilon)\alpha^2]^{|x| - |z_r|}.$$

The right-hand side of (37) is bounded by $\delta$ by (35), since $\alpha < 1$. Therefore, we may apply Lemma 2.7 with (32) and (37) to obtain

$$(38) \qquad D(\mu^{y,S}[k]) \leq (1 + \varepsilon) \sum_{r=1}^{B} D(M\mu^{z_r,S})$$

$$\leq (1 + \varepsilon)\delta\alpha^2 \sum_{r=1}^{B} \sum_{x \in S \cap T(z_r)} [(1 + \varepsilon)\alpha^2]^{|x| - |z_r|}$$

$$= \delta \sum_{x \in S \cap T(y)} [(1 + \varepsilon)\alpha^2]^{|x| - |y|},$$

proving (36).

Applying (36) for the root $\rho$ and $S_n$, we get by (34), taking $n \to \infty$,

$$\limsup_{n \to \infty} D(\mu_1^{\rho,S_n}[k], \ldots, \mu_q^{\rho,S_n}[k]) \leq \delta \lim_{n \to \infty} \sum_{x \in S_n} [(1 + \varepsilon)\alpha^2]^{|x|} = 0,$$

which implies by Lemma 2.8 that

$$\max_{i,j} \lim_{n \to \infty} D_{\mathrm{TV}}(\mu_i^{\rho,S_n}[k], \mu_j^{\rho,S_n}[k]) = 0,$$

as needed.

The proof for $\nu$-robust-solvability is exactly the same. The proof for erasure-robust solvability requires the following modification.

First, as in the proof of Theorem 1.2, we may find $\epsilon^* < 1$ such that if $\epsilon > \epsilon^*$ and $B \leq K$, then for all $y$, which have $B$ children $z_1, \ldots, z_B$ in a cutset $S$, it holds that $D(\mu^{y,S}[\epsilon]) < \delta$ (with $\delta$ as above).

If $S$ is an antichain, let $S'$ denote the set of children of $S$. Note that $S'_n$ is an antichain for all $n$. We prove by induction that for all antichains $S = S_n$ and all $y \in S \cup \mathrm{Ins}(S)$, for $\epsilon > \epsilon^*$,

$$(39) \qquad D(\mu^{y,S'}[\epsilon]) \leq \delta \sum_{x \in S \cap T(y)} [(1 + \varepsilon)\alpha^2]^{|x| - |y|}.$$

The proof is again by induction on $s - |y|$, where $s = \max_{x \in S} |x|$. The only difference is in that for $y \in S$, we use the estimate $D(\mu^{y,S'}[\epsilon]) < \delta$. The remainder of the proof is the same. $\square$

**Acknowledgments.** This research was carried out during a visit to the Computation, Combinatorics and Probability program at the Isaac Newton Institute for Mathematical Sciences in Cambridge, U.K. We thank Yuval Peres for suggesting the 3rd open problem above.

DEPARTMENT OF MATHEMATICS
UPPSALA UNIVERSITY
P.O. BOX 480
SE-751 06 UPPSALA
SWEDEN
E-MAIL: svante.janson@math.uu.se
URL: www.math.uu.se/svante

DEPARTMENT OF STATISTICS
EVANS HALL
UNIVERSITY OF CALIFORNIA, BERKELEY
BERKELEY, CALIFORNIA 94720-1776
USA
E-MAIL: mossel@stat.berkeley.edu
URL: www.stat.berkeley.edu/mossel